\documentclass[12pt,twoside]{article} 
\usepackage[latin1]{inputenc}
\usepackage{graphicx}
\usepackage{enumitem}
\usepackage{latexsym}
\usepackage{amssymb}
\usepackage{mathrsfs}
\usepackage{epsfig}
\usepackage{xcolor}
\usepackage{amsmath}
\usepackage{natbib}
\usepackage{float}
\usepackage{natbib}
\usepackage[singlelinecheck=false]{caption}
\allowdisplaybreaks 

\numberwithin{equation}{section}

\def\beq{\begin{eqnarray*}}
\def\eeq{\end{eqnarray*}}

\begin{document}
\title{\bf Bias Reduction for nonparametric Estimators applied to functional Data Analysis
}

\author{Melanie Birke\\
Department of Mathematics, \\Physics, and Computer Science \\ University of Bayreuth
\and
Tim Greger\\
Department of Mathematics, \\Physics, and Computer Science \\ University of Bayreuth}

\maketitle

\newtheorem{theo}{Theorem}[section]
\newtheorem{lemma}[theo]{Lemma}
\newtheorem{cor}[theo]{Corollary}
\newtheorem{rem}[theo]{Remark}
\newtheorem{prop}[theo]{Proposition}
\newtheorem{defin}[theo]{Definition}
\newtheorem{example}[theo]{Example}
\newtheorem{Assumption}{Assumption}

\begin{abstract}
Compared to nonparametric estimators in the multivariate setting, kernel estimators for
functional data models have a larger order of bias. This is problematic for
constructing confidence regions or statistical tests since the bias might not be negligible. It stems from the fact that one sided kernels are used where already the first moment of the kernel is different from 0. It cannot be cured by assuming the existence of higher order derivatives. In the following, we propose bias corrected estimators based on the idea in \cite{Cheng2018} which still have an appealing structure, but have a bias of smaller order as in multiple regression settings while the variance is of the same order of magnitude as before. In addition we show asymptotic normality of such estimators and derive uniform rates. The performance of the estimator in finite samples is in addition checked in a simulation study.
\end{abstract}



\section{Introduction}
Nonparametric curve estimation is widely used in statistical modeling to estimate complex relationships between variables without assuming a predefined functional form. Unlike parametric models, which impose strict assumptions about the underlying data structure, nonparametric methods allow greater flexibility, making them particularly useful for capturing intricate patterns. However, this flexibility often comes at the cost of increased bias and variance trade-offs, which can compromise the accuracy of the model's predictions.

Bias in nonparametric regression refers to systematic errors that arise when the model's predictions deviate from the true underlying function, often due to inappropriate smoothing or the choice of tuning parameters like bandwidth, \cite{zbMATH00991833}. In practice, bias can lead to underfitting, where the model fails to capture important features of the data. Conversely, reducing bias too aggressively can increase variance, leading to overfitting and poor generalization, \cite{zbMATH02238053}. Therefore, bias reduction is crucial in achieving a well-balanced model that captures the essential patterns in the data without introducing excessive errors.

Several approaches have been proposed to mitigate bias in nonparametric regression, such as optimal bandwidth selection, local polynomial fitting, which mainly improves the behavior of the bias at the boundary, and advanced kernel smoothing techniques, see e.g.~\cite{zbMATH05492743}. These methods aim to strike an effective balance between bias and variance, enhancing the model's predictive power. By understanding and implementing bias reduction techniques, researchers can improve the accuracy and reliability of nonparametric models, which are essential in fields such as economics, machine learning, and biostatistics. 

There is a long history of diverse methods for bias reduction in classical nonparametric regression and related fields. Approaches like local linear estimates in \cite{zbMATH00991833} or the method in \cite{zbMATH01326472} correct the bias only near the boundary. Others reduce the bias on the whole domain. Most of those methods have a two step structure where pilot estimators are modified in a second step to improve the bias. Main contributions in this field are the following.
\cite{zbMATH06036304} extend the method proposed in \cite{Chung2011ALD} for density estimators to the regression framework. The new estimators are based on improving pilot estimators with likelihood methods. They reduce the bias one order of magnitude while preserving the order of variance. Those estimators do not belong to the class of linear smoothers which makes it more difficult to derive asymptotic results.
A two step method, which reduces the bias by multiplying a nonparametrically estimated correction factor to a pilot nonparametric estimator, is given in \cite{zbMATH00786446} for density estimation and discussed in \cite{zbMATH00571679}, \cite{zbMATH01194745} as well as in \cite{zbMATH07064915} for regression. Those methods necessitate the choice of two bandwidths and also have no representation as linear smoothers.
\cite{zbMATH00635765} in a very general setup as well as \cite{Cheng2018} in local polynomial regression use the linear dependence of the bias on some regularisation parameter. Adequately fitting a linear model in a second step reduces the bias. In contrast to the above mentioned methods, the resulting estimators are still linear smoothers which very much simplifies the asymptotic analysis.
In \cite{abramson1982bandwidth} for density estimation and \cite{zbMATH05559558} for binary regression, different bandwidths are chosen for different sample points, which are determined by pilot estimators for the density function. Those approaches are similar to the one proposed in \cite{zbMATH01281011} where, instead of bandwidths, additional weights for the kernels are determined from the estimated density. In contrast, \cite{zbMATH01207185} uses the weighted average of some local linear smoothers to obtain the same order of bias and variance as for the local cubic smoother.
Other methods rely on bootstrap (\cite{zbMATH07031066}), boosting (\cite{DIMARZIO20082483}) or other Monte Carlo methods (\cite{zbMATH07664805}). 
Thorough comparison of different bias optimization techniques, including the one proposed in \cite{zbMATH06036304}, mainly based on Monte Carlo experiments are given in \cite{zbMATH06951373} and \cite{zbMATH06495104}. 
 
In recent years, functional nonparametric estimation has gained more and more impact and many asymptotic results have been derived, see e.g.~\cite{FredericFerraty2006}, \cite{Ferraty2007} or \cite{Ferraty2010}. Due to the use of one sided kernels, the bias, which is of order $h$ when $h$ is the bandwidth, cannot be reduced by assuming higher order differentiablity. None the less, bias reduction has not gained so much attention yet. One reference to mention in this context is \cite{demongeotetal} who construct an analogue of local polynomial estimates in a scalar on function nonparametric setup and derive asymptotic properties as consistency and convergence rates. We also mention \cite{zbMATH07514915}, nevertheless, as model, a Hilbertian on real vector regression is considered there which differs from the one below. There, the bias is reduced in a two step procedure by projecting on a parametric subspace and performing a nonparametric estimation. 

Due to its flexibility and expedient structure, in this paper we generalise the method in \cite{Cheng2018} to more flexible structures of the bias. We also demonstrate, that the method is not limited to regression but can be used for any starting estimate with a certain structure of bias. At the end of Section \ref{se:general} we give some hints on how to choose the regularisation parameters with ideas from optimal design of experiments. As an example, in Section \ref{se:fda}, we apply the method to nonparametric estimation, having functional data as independent variables like scalar on function regression, conditional distribution function or conditonal density. We present asymptotic expansions for bias and variance of the proposed estimator and also show that aymptotic properties of the pilot estimator, like asymptotic distribution or uniform convergence rates, easily transfer to the bias reduced one, which also generalises the results in \cite{Cheng2018}. Finally, we demonstrate the finite sample behavior in a simulation study. 

\section{A general method for bias reduction}
\label{se:general}
Let $\hat m_{h_1,\ldots,h_p}$ be a regularised estimator of some parameter $m$ depending on regularization parameters $h_1,\ldots,h_p$ and let the expectation of the estimator be of the form
\[E[\hat m_{h_1,\ldots,h_p}]= m +\ell_1(h_1)\beta_1+\ldots+\ell_p(h_p)\beta_p+o(\max_{i=1,\ldots,p}\ell_i(h_i)).\]
\subsection{Definition and first properties}
If we take $\hat m_{h_1,\ldots,h_p}$ as an approximation of its expectation and take $B$ different choices $\mathbf h_i=(h_{i1},\ldots,h_{ip})$, $i=1,\ldots,B$ of regularization vectors, we approximatly get a linear model
\[\hat m_{\mathbf h_i}\approx \beta_0+\ell_1(h_{i1})\beta_1+\ldots+\ell_p(h_{ip})\beta_p+\delta_i\]
for $\beta_0=m$ and some error term $\delta_i$. With $V^T=(\hat m_{\mathbf h_1},\ldots,\hat m_{\mathbf h_B})$, $e_1^T=(1,0,\ldots,0)\in \mathbb R^{1\times p+1}$ and

\begin{align*}
H&=\begin{pmatrix}
1&\ell_1(h_{11})&\cdots&\ell_p(h_{1p})\\
\vdots&\vdots&&\vdots\\
1&\ell_1(h_{B1})&\cdots&\ell_p(h_{Bp})
\end{pmatrix}
\end{align*}
we get a further estimator of $m$ as
\begin{align}
\label{eq:biasredest}
\hat m_{B}&=e_1^T(H^TH)^{-1}H^TV.
\end{align}

With this we have
\begin{align}\label{properties}
e_1^T(H^TH)^{-1}H^T1_{B}&=1\nonumber\\
e_1^T(H^TH)^{-1}H^T(\ell_j(h_{1j}),\ldots,\ell_j(h_{Bj}))^T&=0,\ j=1,\ldots,p.\nonumber\\
\end{align}
This can be easily seen, since $H^T1_B$ is the first column of the matrix $H^TH$, such that $e_1^T(H^TH)^{-1}H^T1_{B}$ is the $(1,1)$th entry of the identity matrix. In the same way we see, that $H^T(\ell_j(h_{1j}),\ldots,\ell_j(h_{Bj}))^T$ is the $(j+1)$th column of $H^TH$ and therefore $$e_1^T(H^TH)^{-1}H^T(\ell_j(h_{1j}),\ldots,\ell_j(h_{Bj}))^T$$ is the $(1,j+1)$th entry of the identity matrix, $j=1,\ldots,p$ which is $0$. From (\ref{properties}) the next theorem can easily be deduced.

\begin{theo}
\label{th:biasred_general}
Under the above assumptions there is
\[E[\hat m_{B}]-m=o(\max_{i=1,\ldots,p}\ell_i(h_i)).\]
\end{theo}
It is also possible to transfer other asymptotic properties as order of variance, convergence rates or central limit theorems from the initial estimator to the bias reduced one. But this depends on the special structure of the initial estimator and the underlying model. Therefore we discuss this exemplarily for nonparametric functional data models in the next section.

If we take $m=r(\cdot):\mathbb R\to\mathbb R$ and consider the model
\[Y_i=r(X_i)+\varepsilon_i,\ i=1,\ldots,n\]
with $r$ twice continuously differentiable, $(Y_i,X_i)$ independent identically distributed random variables with values in $\mathbb R^2$, where $X_i$ has twice continuously differentiable density $f$, $\varepsilon_i$ error terms with $E[\varepsilon_i]=0$ and $Var(\varepsilon_i)<\infty$, $i=1,\ldots,n$, we obtain the estimator proposed in \cite{Cheng2018} by choosing the local linear estimator $\hat m=\hat r_{h_1}(x)$ in the fixed argument $x\in\mathbb R$ as initial estimator and further $p=1$, $\ell_1(h_1)=h_1^2$ for the bias reduction.

\subsection{Choosing the Sequence of Regularization Parameters}
As we will see later on, the method above reduces the order of the bias, but might increase the variance compared to the original estimator due to a multiplicative constant. Fortunately, the order of the variance stays the same but one might be interested in choosing the sequence of regularisation parameters in a way that the variance is smallest. Here we can make use of results from optimal design of experiments, which is well established for linear models, see e.g. \cite{zbMATH05022305}. From this theory we get for fixed $B$ a design consisting of points $k_1,\ldots,k_D$, $1\leq D\leq B$ where to observe and weights $p_1,\ldots,p_D$ giving the proportion of observations in the respective points. All this relies on the assumption that repeated measurements are possible and result in different outcomes. This is not the case in our situation, since, given a sample, for one explicite regularisation parameter we always get the same value for the estimator. Guided by the theoretical result, one possibility is to choose the $p_jB$ observations closely around the points $k_j$, $j=1,\ldots,D$ instead of exactly in those points. What 'closely' means, very much depends on the size of the regularisation parameters and the model itself. We will deeper discuss this in the particular examples in the simulation section.

\section{Application to Nonparametric Functional Data Analysis}
\label{se:fda}
Let $(Y_i,\mathcal{X}_i)$, $i=1,\ldots,n$ be a sample of i.i.d. pairs of random variables. The random variables $(Y_i)_{i=1,\ldots,n}$ are real valued and $(\mathcal{X}_i)_{i=1,\ldots,n}$ are random elements with values in a functional space $\mathcal{E}$. We take $\mathcal{E}$ to be a separable Banach space endowed with a norm $||\cdot||$. One possibility is to consider the regression model
\begin{align*}
Y_k=r(\mathcal{X}_k)+\varepsilon_k,~~~k=1,...,n
\end{align*}
where $r:\mathcal E\to\mathbb R$ is the regression function while $\varepsilon_1,\ldots,\varepsilon_n$ are independent identically distributed with $E[\varepsilon_k|\mathcal{X}_k]=0$ as well as $E[\varepsilon_k^2|\mathcal{X}_k]=\sigma_{\varepsilon}^2(\mathcal{X}_k)<\infty$. Other quantities we consider are conditional distribution functions and conditional densities of $Y_1$ given $\mathcal X_1$.

All estimators we consider are of a general form
\begin{align*}
\widehat{m}_\Phi(\chi)=\frac{\sum_{k=1}^n \Phi_n(Y_k)K(h^{-1}||\mathcal{X}_k-\chi||)}{\sum_{k=1}^n K(h^{-1}||\mathcal{X}_k-\chi||)}.
\end{align*}
Here $K:[0,1]\to\mathbb R^+$ is a onesided kernel as in \cite{FredericFerraty2006}, $h=h_n>0$ is a bandwidth and $\chi\in \mathcal{E}$ is a fixed function. Further assumptions on these quantities will be stated later on. The function $\Phi_n$ may depend on $n$ but not on $h$.
As we continue, we will look at 
\begin{enumerate}
\item $\Phi_n(Y_k)=Y_k$ for estimating the regression function.
\item $\Phi_n(Y_k)=I_{\lbrace Y_k\leq y\rbrace}$ for estimating the conditional distribution function.
\item $\Phi_n(Y_k)=\frac1bK_0(b^{-1}(y-Y_k))$, $y\in \mathbb{R}$, a symmetric kernel $K_0:[-1,1]\to\mathbb R_0^+$ and a bandwidth $b>0$ depending on $n$ which goes to zero as $n$ goes to infinity, for estimating the conditional density.
\end{enumerate}
As stated in \cite{FredericFerraty2006}, these three cases give us nonparametric estimators for the regression function $r(\chi)$, the conditional distribution function $F_Y^{\mathcal{X}}(\chi,y)=P(Y\leq y|\mathcal{X}=\chi)=E[I_{\lbrace Y\leq y\rbrace}|\mathcal{X}=\chi]$ and the conditional density function $f_Y^{\mathcal{X}}(\chi,y)=\frac{\partial}{\partial y}F_Y^{\mathcal{X}}(\chi,y)$, which will be denoted as $\widehat{r}(\chi)$, $\widehat{F}_Y^{\mathcal{X}}(\chi,y)$ and $\widehat{f}_Y^{\mathcal{X}}(\chi,y)$. But even other choices of $\Phi_n$ are covered by the following results as long as adequate assumptions are fulfilled. For a unified notation assume, that $E[\Phi_n(Y_1)|\mathcal X_1=\chi]$ exists and can be wirtten as
\[E[\Phi_n(Y_1)|\mathcal X_1=\chi]=m_\Phi(\chi)+R_n(\chi,\Phi)\]
for $n$ to infinity. For the three examples from above, this is
\begin{enumerate}
\item $m_\Phi=r$ for $\Phi_n(Y_k)=Y_k$.
\item $m_\Phi=F_Y^{\mathcal X}(\cdot,y)$ for $\Phi_n(Y_k)=I_{\lbrace Y_k\leq y\rbrace}$.
\item $m_\Phi=f_Y^{\mathcal X}(\cdot,y)$ for $\Phi_n(Y_k)=\frac1bK_0(b^{-1}(y-Y_k))$, $y\in \mathbb{R}$.
\end{enumerate}
We define in general
\[\varphi_{\Phi,\chi}(s)=E[m_\Phi(\mathcal X)-m_\Phi(\chi)|~||\mathcal{X}-\chi||=s]\]
and for the explicit examples
\begin{align*}
\varphi_{reg,\chi}(s)=E[(r(\mathcal{X})-r(\chi))~\vert~||\mathcal{X}-\chi||=s],
\end{align*}
\begin{align*}
\varphi_{cdf,\chi}(s)=E[(F_Y^{\mathcal{X}}(\mathcal{X},y)-F_Y^{\mathcal{X}}(\chi,y))~|~||\mathcal{X}-\chi||=s],
\end{align*}
\begin{align*}
\varphi_{dens,\chi}(s)=E[(f_Y^{\mathcal{X}}(\mathcal{X},y)-f_Y^{\mathcal{X}}(\chi,y))|~||\mathcal{X}-\chi||=s].
\end{align*}
Furthermore, the small ball probability is set as
\begin{align*}
L_\chi(t)=P(||\mathcal{X}-\chi||\leq t),
\end{align*}
and there is
\begin{align*}
\tau_{h,\chi}(s)=\frac{L_\chi(hs)}{L_\chi(h)}=P(||\mathcal{X}-\chi||\leq hs~\vert~||\mathcal{X}-\chi||\leq h)~~~\forall s\in [0,1].
\end{align*}
Obviously $L_\chi(t)$ is the distribution function of the random variable $||\mathcal{X}-\chi||$. We also have the following assumptions.
\begin{itemize}
\item[$A_0$:] $m_\Phi$ is continuous in a neighbourhood of $\chi\in\mathcal E$.
\item[$A_1'$:] $\varphi_{\Phi,\chi}'(0)$ exists.
\item[$A_2$:] $\lim_{n\rightarrow \infty}h=0$, $L_\chi(0)=0$ and $\lim_{n\rightarrow \infty} nL_\chi(h)=\infty$.
\item[$A_3$:] $\forall s\in [0,1]: \tau_{h,\chi}(s)\rightarrow \tau_{0,\chi}(s)$ for $h\rightarrow 0$.
\item[$A_4$:] $K$ is supported on $[0,1]$ with $K(1)>0$ and has a a derivative on $[0,1)$ with $-\infty<C_1<K'(t)<C_2<0$ for some constants $C_1,~C_2$ and $t\in[0,1)$.
\item[$A_5$:] $K_0:[-1,1]\to\mathbb R_0^+$ is a second order Lipschitz continuous symmetric kernel with $\int_{-1}^1K_0^2(t)dt<\infty$.

\end{itemize}
Remark that $\varphi_{\Phi,\chi}(0)=0$. The measurable mappings $s\to\tau_{h,\chi}(s)$ and $s\to\tau_{0,\chi}(s)$ are nondecreasing in $s$. Explicit forms of $\tau_{0,\chi}$ for different types of functional data are given in \cite{Ferraty2007}.
\\In addition define
\begin{align*}
M_{0,\chi}=K(1)-\int_0^1 (sK(s))'\tau_{0,\chi}(s)ds
\end{align*}
and
\begin{align*}
M_{1,\chi}=K(1)-\int_0^1 K'(s)\tau_{0,\chi}(s)ds.
\end{align*}
The following result generalises the results for the order of magnitude of the bias in \cite{Ferraty2007} to a general estimator  $\hat m_\Phi(\chi)$, which includes e.g.~estimators of $F_Y^{\mathcal{X}}(\chi,y)$ or $f_Y^{\mathcal{X}}(\chi,y)$.
\begin{theo}\label{Ferraty07T1}
Under assumptions $A_0$, $A_1'$, $A_2-A_4$, if $\chi\in \mathcal E$ is fixed, we have
\begin{align*}
Bias(\widehat{m}_\Phi(\chi))=\varphi_{\Phi,\chi}'(0)\frac{M_{0,\chi}}{M_{1,\chi}}h+R_n(\chi,\Phi)+O((nL_\chi(h))^{-1})+o(h).
\end{align*}
\end{theo}
It is easy to verify that for $\Phi_n=id$ as well as $\Phi_n=I\{\cdot\leq y\}$, there is $R_n(\chi,\Phi)=0$ while under condition $A_5$, if $\frac{\partial^2}{\partial y^2}{f_Y^{\mathcal X}}(\chi,y)$ exists and is continuous, for $\Phi_n=\frac1bK_0(b^{-1}(\cdot-y))$ there is $R_n(\chi,\Phi)=b^2\frac12\int_{-1}^1v^2K_0(v)dv\frac{\partial^2}{\partial y^2}{f_Y^{\mathcal X}}(\chi,y)+o(b^2)$. This directly results in the following corollary.

\begin{cor}
Under assumptions $A_0$, $A_1'$, $A_2-A_4$, if $\chi\in \mathcal E$ and $y\in\mathbb R$ are fixed, we have
\begin{enumerate}
\item for regression estimation
\begin{align*}
Bias(\widehat{r}(\chi))=\varphi_{reg,\chi}'(0)\frac{M_{0,\chi}}{M_{1,\chi}}h+O((nL_\chi(h))^{-1})+o(h).
\end{align*}
\item for estimating the conditional distribution function in $y\in\mathbb R$ fixed
\begin{align*}
Bias(\widehat{F}_Y^{\mathcal{X}}(\chi,y))=\varphi_{cdf,\chi}'(0)\frac{M_{0,\chi}}{M_{1,\chi}}h+O((nL_\chi(h))^{-1})+o(h).
\end{align*}
\item for estimating the conditional density function, if additionally $A_5$ is true and $\frac{\partial^2}{\partial y^2}{f_Y^{\mathcal X}}(\chi,y)$ exists and is continuous
\begin{align*}
Bias(\widehat{f}_Y^{\mathcal{X}}(\chi,y))=&\varphi_{pdf,\chi}'(0)\frac{M_{0,\chi}}{M_{1,\chi}}h+\frac12\frac{\partial^2}{\partial y^2}{f_Y^{\mathcal X}}(\chi,y)\int_{-1}^1v^2K_0(v)dv~b^2\\&+O((nL_{\chi}(h))^{-1})+o(h)+o(b^2).
\end{align*}
\end{enumerate}
\end{cor}

We will now construct bias-reduced versions of these estimators with the method proposed in Section \ref{se:general}.
Choose a set of bandwidths $\mathbf{h}_1,...,\mathbf{h}_B$ with $\mathbf{h}_i=h_i$ such that $h_1<\ldots<h_B$ if the reduction shall only be in $h$, or $\mathbf{h}_i=(h_i,b_i)$ with $h_i$ as above and $b_1<\ldots<b_B$, if the reduction is desired to be both in $h$ and $b$. 
The bias reduced estimator $\widehat m_{\Phi,B}(\chi)$ is then computed as in (\ref{eq:biasredest}) with 
\[H_1=\begin{pmatrix}
1&h_1\\
\vdots&\vdots\\
1&h_B
\end{pmatrix}\]
for the reduction only in $h$, which can be applied for estimating all functions $m_\Phi$, and
\[H_2=\begin{pmatrix}
1&h_1&b_1\\
\vdots&\vdots&\vdots\\
1&h_B&b_B
\end{pmatrix}\]
for reducing in $h$ and $b$ which is only reasonable for estimators like the conditional density estimator. For both matrices, the estimators can be written as
\[\widehat m_{\Phi,B}(\chi)=\sum_{i=1}^Bg_{i,j}\widehat m_{\Phi,\mathbf h_i}(\chi),~j=1,2\]
with weights $g_{i,j}$, $i=1,\ldots,B$, $j=1,2$, which are given by
\[g_{i,1}=\frac{\sum_{k=1}^B h_k^2-h_i\sum_{k=1}^B h_k}{B\sum_{k=1}^B h_k^2-\left(\sum_{k=1}^B h_k\right)^2}
\]
for $H_1$ but having a more complicated form for $H_2$, which we do not state here explicitly. Nevertheless, $g_{i,2}$ can easily be deduced by matrix calculations using $H_2$.
In the sequel the bias reduced estimators for the explicit examples will be denoted as $\widehat r_B(\chi)$, $\widehat{F}_{Y,B}^{\mathcal{X}}(\chi,y)$ and $\widehat{f}_{Y,B}^{\mathcal{X}}(\chi,y)$ for the reduction only in $h$. Furthermore, we write $\widehat{f}_{Y,B,2}^{\mathcal{X}}(\chi,y)$ for the conditional density estimator with reduction in $h$ and $b$.

\subsection{Bias and Variance}
In this section we give asymptotic results concerning bias and variance for estimators  of the form $\widehat m_{\Phi,B}$. For the higher order expansion of the bias we assume
\begin{itemize}
\item [$A_1$:] $\varphi_{\Phi,\chi}''(0)$ exists.
\item [$A_6$:] $h_i=C_ih_0$ for $i=1,...,B$ with $1< C_1< ...< C_B$ and some $h_0>0$ as well as $b_i=C'_ib_0$ $\forall i=1,...,B$ with $1< C'_1< ...< C'_B$ and some $b_0>0$.
\end{itemize}
Under assumption $A_6$ the weights $g_{i,j}$, $j=1,2$ are independent of the bandwidths $h_0$ and $b_0$, e.g. for $H_1$
\[g_{i,1}=\frac{\sum_{k=1}^B C_k^2-C_i\sum_{k=1}^B C_k}{B\sum_{k=1}^B C_k^2-\left(\sum_{k=1}^B C_k\right)^2}.\]
Furthermore, we define
\[M_{3,\chi}=\int_0^1(s^2K'(s)+2sK(s))\tau_{0,\chi}(s)ds.\]
The following theorem states results for the bias with the reduction only in the regularisation parameter $h$.
\begin{theo}\label{TheoBias}
Under the assumptions $A_0-A_4$ and $A_6$, we have for $\chi\in\mathcal E$ and $y\in\mathbb R$ fixed

\begin{align*}
Bias(\widehat{m}_{\Phi,B}(\chi))&=\frac12\sum_{i=1}^Bg_{i,1}C_i^2\varphi_{\Phi,\chi}''(0)\frac{M_{3,\chi}}{M_{1,\chi}}h_0^2+R_n(\chi,\Phi)+O((nL_{\chi}(h))^{-1})+o(h_0^2)\\
&=S(\chi,\Phi)h_0^2+R_n(\chi,\Phi)+O((nL_{\chi}(h))^{-1})+o(h_0^2)
\end{align*}
with $S(\chi,\Phi)=\frac12\sum_{i=1}^Bg_{i,1}C_i^2\varphi_{\Phi,\chi}''(0)\frac{M_{3,\chi}}{M_{1,\chi}}$.
\end{theo}
For the three examples from above we again state the results as corollary.
\begin{cor}
\label{cor:biasred}
Under assumptions $A_0$, $A_1'$, $A_2-A_4$ and $A_6$, if $\chi\in \mathcal E$ is fixed, we have
\begin{enumerate}
\item for regression estimation
\begin{align*}
Bias(\widehat{r}_{B}(\chi))&=\frac12\sum_{i=1}^Bg_{i,1}C_i^2\varphi_{reg,\chi}''(0)\frac{M_{3,\chi}}{M_{1,\chi}}h_0^2+O((nL_{\chi}(h))^{-1})+o(h_0^2).
\end{align*}
\item for the estimator of the conditional distribution function in $y\in\mathbb R$
\begin{align*}
Bias(\widehat{F}_{Y,B}^{\mathcal{X}}(\chi,y))=\frac12\sum_{i=1}^Bg_{i,1}C_i^2\varphi_{cdf,\chi}''(0)\frac{M_{3,\chi}}{M_{1,\chi}}h_0^2+O((nL_{\chi}(h))^{-1})+o(h_0^2).
\end{align*}

\item for estimating the conditional density function, if additionally $A_5$ is true and $\frac{\partial^2}{\partial y^2}{f_Y^{\mathcal X}}(\chi,y)$ exists 
\begin{align*}
Bias(\widehat{f}_{Y,B}^{\mathcal{X}}(\chi,y))=&\frac12\sum_{i=1}^Bg_{i,1}C_i^2\varphi_{dens,\chi}^{''}(0)\frac{M_{3,\chi}}{M_{1,\chi}}~h_0^2
\\&+\frac12\frac{\partial^2}{\partial y^2}{f_Y^{\mathcal{X}}}(\chi,y)\int_{\mathbb{R}}v^2K_0(v)dv~b^2\\&+O((nL_{\chi}(h))^{-1})+o(h_0^2)+o(b^2).
\end{align*}
\end{enumerate}

\end{cor}

\begin{rem}
\label{rem:Bchoice}
If, under the conditions of Theorem \ref{TheoBias} the intervall of the bandwidth taken is fixed to $(h_0,h^u]$ for some $h^u>h_0$ and increase the number $B$ of equally spaced bandwidth $h_i=h_0\left(1+i\frac{h^u-h_0}{h_0B}\right)$, $i=1,\ldots,B$, we easily see, that
\[\sum_{i=1}^Bg_{i,1}C_i^2\to -\left(1+\frac{h^u-h_0}{h_0}+\frac16\left(\frac{h^u-h_0}{h_0}\right)^2\right).\]
This means, that increasing the number of bandwidths used for the bias reduction can only marginally improve the reduction. This is, what we also observe in Tables \ref{Tab:8} and \ref{Tab:9}.
\end{rem}

Regarding the conditional density function it is possible to apply the bias reduction on both bandwidths $h$ and $b$, using the result from chapter \ref{se:general}.

\begin{lemma}\label{lemma}
Under assumptions $A_0-A_6$, if $\frac{\partial^4}{\partial y^4}f_Y^{\mathcal{X}}(\chi,y) $ exists
\begin{align*}
Bias(\widehat{f}_{Y,B,2}^{\mathcal{X}}(\chi,y))&=h_0^2\frac12\varphi_{dens,\chi}^{''}(0)\frac{M_{3,\chi}}{M_{1,\chi}}\sum_{i=1}^Bg_{i,2}C_i^2+o(h_0^2)
\\&+b_0^4\frac{1}{24}f_Y^{(4),\chi}(y)\int_{\mathbb{R}}v^4K_0(v)dv\sum_{i=1}^Bg_{i,2}C_i^{'4}+O((nL_{\chi}(h))^{-1})+o(h_0^2)+o(b_0^4),
\end{align*}
\end{lemma}

\begin{rem}
The main difference between $\widehat{f}_{Y,B,2}^{\mathcal{X}}(\chi,y)$ and $\widehat{f}_{Y,B}^{\mathcal{X}}(\chi,y)$ is the choice of the weights $g_{i,2}$ and $g_{i,1}$. The significant properties of these weights will be crucial only regarding the bias. Computations for variance, asymptotic normality and uniform rates can be done in exactly the same way. Therefore we just highlight the double bias-reduction of $\widehat{f}_{Y,B,2}^{\mathcal{X}}(\chi,y)$ and refer for variance, asymptotic normality and uniform rates to $\widehat{f}_{Y,B}^{\mathcal{X}}(\chi,y)$.\\
In this situation the reduction of the bias with respect to $b$ can be obtained with simpler methods by just choosing a kernel $K_0$ of fourth order instead of second order. This reduces the computational effort and even results in a simpler structure of the bias in $b$ without the constant $\sum_{i=1}^Bg_{i,2}C_i^{'4}$. Nevertheless the simultaneous bias reduction in two or more parameters might be of higher interest in other settings. 
\end{rem}

In the following result we are going to state the asymptotic variances. To this end we need the additional assumption:
\begin{itemize}
\item[$A_7$:] $E[\Phi_n^2(Y_1)]<\infty$ and $D_{\Phi,n}(\chi)=Var(\Phi_n(Y_1)|\mathcal X_1=\chi)$ is continuous in $\chi$.
\end{itemize}
\begin{theo}\label{TheoVar}
Under the assumptions $A_0-A_4$, $A_6$ and $A_7$ there is for fixed $\chi\in\mathcal E$, $y\in\mathbb R$ and
\begin{align*}
M_{2,i_1,i_2,\chi}=K(1)K(C_{i_2}/C_{i_1})-\int_0^1K(s)K'\Bigl(s\frac{C_{i_2}}{C_{i_1}}\Bigr)\tau_{0,\chi}(s)ds
\end{align*}
\begin{align*}
Var(\widehat{m}_{\Phi,B}(\chi))=\frac{1}{nL_\chi(h_0)}\frac{D_{\Phi,n}(\chi)}{M_{1,\chi}^2}\sum_{i_1,i_2=1}^Bg_{i_1,j}g_{i_2,j}\tau_{0,\chi}(1/C_{i_1})M_{2,i_1,i_2,\chi}+o\Bigl(\frac{D_{\Phi,n}(\chi)}{nL_{\chi}(h_0)}\Bigr),~j=1,2.
\end{align*}
\end{theo}

Remark that these variances are of the same order as the variances of the classical estimators, see e.g. \cite{Ferraty2007}.
With $D_{\Phi,n}(\chi)=\sigma^2(\chi)$, $D_{\Phi,n}(\chi)=F_Y^{\mathcal X}(\chi,y)(1-F_Y^{\mathcal X}(\chi,y))$ respectively $D_{\Phi,n}(\chi)=\frac1b\int_{-1}^1 K_0^2(v)dvf_Y^{\mathcal X}(\chi,y)+o(1/b)$, this yields the corresponding results for $\hat r_B$, ${\hat F}_Y^{\mathcal X}(\chi,y)$ respectively ${\hat f}_Y^{\mathcal X}(\chi,y)$.

\subsection{Asymptotic Normality and Uniform Rates}
Having derived the pointwise asymptotic decompositions of bias and variance we are now in a position to derive  further asymptotic  results.
\subsection{Asymptotic Distribution}
We only state the general result here. The explicit cases for regression, conditional distribution and conditional density function can easily be deduced from the results in the last section.
\begin{theo}\label{TheoAsymNorm}
Under the same assumptions as in Theorem \ref{TheoVar} and if $E[|\Phi_n(Y_1)|^4]<\infty$ for all $n\in\mathbb{N}$, $\frac{nL_\chi(h_0)}{D_{\Phi,n}(\chi)}h_0^4=O(1)$ and $\frac{nL_\chi(h_0)}{D_{\Phi,n}(\chi)}R_n^2(\chi,\Phi)=O(1)$ for $n \to\infty$, we have

\begin{align*}
\frac{\sqrt{nL_\chi(h_0)}}{\sqrt{D_{\Phi,n}(\chi)}}\Bigl(\widehat{m}_{\Phi,B} (\chi)-m_\Phi(\chi)-S(\chi,\Phi)h_0^2-R_n(\chi,\Phi)\Bigr)\overset{\mathcal{D}}{\rightarrow}\mathcal{N}(0,\Gamma)
\end{align*}
with
\begin{align*}
\Gamma=\frac{1}{M_{1,\chi}^2}\sum_{i_1,i_2=1}^Bg_{i_1,j}g_{i_2,j}\tau_{0,\chi}(1/C_{i_1})M_{2,i_1,i_2,\chi},~j=1,2,
\end{align*}
$S(\chi,\Phi)$, $R_n(\chi,\Phi)$ as in Theorem \ref{TheoBias} and $D_{\Phi,n}(\chi)$ as in Theorem \ref{TheoVar}.
\end{theo}
In the case of conditional density estimation it is possible to choose $b=o(h_0)$ or $b=o(h_0^2)$ to get rid of the second bias term. This is of particular advantage if the above results are used for e.g.~constructing pointwise asymptotic confidence intervals.
\subsection{Uniform Rates}
We now turn to the uniform rates of the bias reduced estimators. In a first step, we state the general result for uniform rates in $\chi$ while the possible second argument $y$ is kept fixed. Since the uniformity in $y$ only makes sense for the conditional distribution and conditional density function, we state those results separately. We need all assumptions stated above uniformly in $\chi$ (and $y$) and some additional ones for deriving the uniform almost complete convergence. The results generalise the ones in \cite{Ferraty2010} to the bias reduced ones. Assume $S_{\mathcal{E}}$ to be some subset of $\mathcal{E}$. 
\begin{itemize}
\item[$B_0$:] Assumptions $A_0-A_6$ are fulfilled uniformly in  $\chi\in S_{\mathcal E}$.

\item[$B_1$:]  There are constants $0<\kappa_1\leq\kappa_2<\infty$ and a function $\ell:\mathbb R_0^+\to\mathbb R_0^+$ with existing bounded derivative on some subset $[0,\eta_0)$.
\begin{enumerate}
\item $\kappa_1\ell(h)<L_\chi(h)<\kappa_2 \ell(h)$ for all $h>0$.
\item $\nu_0(s):=\lim_{h\to0}\frac{\ell(hs)}{\ell(h)}$ exists.
\end{enumerate}

\item[$B_2$:] $\forall m\geq 2,~E[|\Phi_n(Y_1)|^m|\mathcal{X}_1=\chi]<\delta_m(\chi)<C<\infty$ with $\delta_m(\cdot)$ continuous on $S_{\mathcal{E}}$.

\item[$B_3$:] Kolmogorov's $\varepsilon$-entropy $\Psi_{S_{\mathcal E}}(\varepsilon)=\log(N_{\varepsilon}(S_{\mathcal E}))$ with $N_{\varepsilon}(S_{\mathcal E})$ the covering number of $S_{\mathcal E}$ fulfills 
\begin{align*}
\sum_{n=1}^{\infty}&\exp\Bigl((1-\beta)\Psi_{S_{\mathcal{E}}}(\log n/n)\Bigr)< \infty\\
(\log n)^2/nd_{\Phi,n}\ell(h)&<\Psi_{S_{\mathcal{E}}}(\log n/n)<nd_{\Phi,n}\ell(h)/\log n
\end{align*}
for some $\beta>1$, if $n$ is large enough, where there are constants $0<c_1<c_2<\infty$ such that $\frac{c_1}{d_{\Phi,n}}<D_{\Phi,n}(\chi)<\frac{c_2}{d_{\Phi,n}}$ for all $\chi\in S_{\mathcal E}$ and $nd_{\Phi,n}\ell(h)\to\infty$.
\end{itemize}
Because of Assumption $B_1$ a), it holds that
\begin{align}\label{EWAbschatzung}
\forall \chi\in S_{\mathcal{E}},\exists 0<\kappa_1<\kappa_2<\infty,~~\kappa_1\ell(h)<E[K(h^{-1}||\mathcal{X}-\chi||)]<\kappa_2\ell(h).
\end{align}
In addition, we obtain from $B_1$ b) that for all $s\in[0,1]$
\[\sup_{\chi\in S_{\mathcal E}}\tau_{0,\chi}(s)\leq \frac{\kappa_2}{\kappa_1}\nu_0(s).\]
In condition $B_3$, $d_{\Phi,n}=1$ for $\Phi_n=id$ as well as $\Phi_n=I\{\cdot\leq y\}$ while $d_{\Phi,n}=\frac1b$ if $\Phi_n=\frac1bK_0\left(b^{-1}(\cdot-y)\right)$. 
The following theorem states the uniform rates for the general estimator $\hat m_{\Phi,B}$ in the functional argument. The explicit cases are discussed in the sequel, also considering the uniform rates in both, the functional and the real argument.
\begin{theo}\label{GKT}
Under $B_0$-$B_3$, if $\sup_{\chi\in S_{\mathcal E}}|R_n(\chi,\Phi)|<\infty$ we have
\begin{align*}
\underset{\chi\in S_{\mathcal{E}}}{\sup}|\widehat{m}_{\Phi,B}(\chi)-m_\Phi(\chi)|=O(h_0^2)+\sup_{\chi\in S_{\mathcal E}}|R_n(\chi,\Phi)|+O_{a.co.}\left(\frac{\Psi_{S_{\mathcal E}}({\log n/n})}{nd_{\Phi,n}\ell(h_0)}\right)^{1/2}.
\end{align*}
\end{theo}
In case of $\Phi_n=id$ or $\Phi_n=I\{\cdot\leq y\}$ for some $y\in\mathbb R$ fixed, $R_n(\chi,\Phi)=0$ as already stated above and a direct consequence of this result is
\begin{align*}
\underset{\chi \in S_{\mathcal{E}}}{\sup}|\widehat{r}_{B}(\chi)-r(\chi)|=O(h_0^2)+O_{a.co.}\left(\frac{\Psi_{S_{\mathcal E}}({\log n/n})}{n\ell(h_0)}\right)^{1/2}=\underset{\chi \in S_{\mathcal{E}}}{\sup}|\widehat{F}_{Y,B}^{\mathcal{X}}(\chi,y)-F_Y^{\mathcal{X}}(\chi,y)|.
\end{align*}
Note, that for the conditional distribution function condition $B_2$ is automatically fulfilled. For the regression case, Condition $B_2$ is the existence, uniform boundedness and continuity of all absolute conditional moments of $Y_1$.
For the conditional density estimator, if $\sup_{\chi\in S_{\mathcal E}}\Big|\frac{\partial^2}{\partial y^2}f_Y^{\mathcal X}(\chi,y)\Big|<\infty$ and $K_0$ is bounded, there is
\begin{align*}
\underset{\chi \in S_{\mathcal{E}}}{\sup}|\widehat{f}_{Y,B}^{\mathcal{X}}(\chi,y)-f_Y^{\mathcal{X}}(\chi,y)|=O(h_0^2)+O(b^2)+O_{a.co.}\left(\frac{\Psi_{S_{\mathcal E}}({\log n/n})}{nb\ell(h_0)}\right)^{1/2},
\end{align*}
which states the uniform rates in the functional argument. \\
\begin{table}[H]
\begin{tabular}{|c|cc|cc|cc|}
\hline 
\multicolumn{1}{|c|}{}&\multicolumn{2}{|c|}{Squared Bias}&\multicolumn{2}{|c|}{MSE}&\multicolumn{2}{|c|}{Variance}\\~&
$\widehat{r}$& $\widehat{r}_B$&$\widehat{r}$&$\widehat{r}_B$&$\widehat{r}$&$\widehat{r}_B$\\\hline 
$n=100$&11.09873&9.330857&11.20392&9.439501&0.1051811&0.1086443\\\hline
$n=200$&10,52505&7.426841& 10.58222&7.486948& 0.05717157&0.06010725\\\hline
$n=500$&8.999492&2.545537&9.025379&2.592098&0.02588676&0.04656156
\\\hline
\end{tabular}
\caption{Squared Bias, MSE and Variance of $\hat r$ with bandwidth $h=9.28\cdot n^{-1/3}$ given in the left column and of $\hat r_B$ using $B=21$ bandwidths centered around $h$ having distance of 0.01 to each other.}
\label{Tab:1}
~
\medskip
\\\begin{tabular}{|c|cc|cc|cc|}
\hline 
\multicolumn{1}{|c|}{}&\multicolumn{2}{|c|}{Squared Bias}&\multicolumn{2}{|c|}{MSE}&\multicolumn{2}{|c|}{Variance}\\~&
$\widehat{r}$& $\widehat{r}_B$&$\widehat{r}$&$\widehat{r}_B$&$\widehat{r}$&$\widehat{r}_B$\\\hline
$n=100$&10.6491&8.222701&10.76833&8.342714&0.1192261&0.1200131\\\hline
$n=200$&9.98805&5.682518&10.04799&5.75525&0.05994118&0.07273146\\\hline
$n=500$&7.609054&0.1898303&7.632157&0.2555752&0.02310295&0.06574486
\\\hline
\end{tabular}
\caption{Squared Bias, MSE and Variance of $\hat r$ with bandwidth $h=8.12\cdot n^{-1/3}$ given in the left column and of $\hat r_B$ using $B=21$ bandwidths centered around $h$ having distance of 0.01 to each other.}
\label{Tab:2}
~
\medskip
\\\begin{tabular}{|c|cc|cc|cc|}
\hline 
\multicolumn{1}{|c|}{}&\multicolumn{2}{|c|}{Squared Bias}&\multicolumn{2}{|c|}{MSE}&\multicolumn{2}{|c|}{Variance}\\~&
$\widehat{r}$& $\widehat{r}_B$&$\widehat{r}$&$\widehat{r}_B$&$\widehat{r}$&$\widehat{r}_B$\\\hline 
$n=100$&10.2908&6.757813&10.39609&6.882673&0.1052884&0.1248599\\\hline
$n=200$&9.167143&2.992993&9.221726&3.080701& 0.05458331&0.08770772\\\hline
$n=500$&5.33639&1.553648&5.361183&1.665007&0.02479364&0.1113599
\\\hline
\end{tabular}
\caption{Squared Bias, MSE and Variance of $\hat r$ with bandwidth $h=6.96\cdot n^{-1/3}$ given in the left column and of $\hat r_B$ using $B=21$ bandwidths centered around $h$ having distance of 0.01 to each other.}
\label{Tab:3}
\end{table}

To get the uniform rates in both the functional and nonfunctional argument, let $S_{\mathbb{R}}\subset \mathbb{R}$ be a compact subset. In our examples, this only makes sense for the conditional distribution function as well as the conditional density functions. There might be other situations like e.g. for the hazard function, see e.g. \cite{Ferraty2010} and conditons on a general $\Phi_n$ depending on $y$ can easily be derived along the lines, but would additionally complicate the notation.
\begin{itemize}
\item[$B_4$:] $F_Y^{\mathcal X}(\chi,y)$ is Hoelder continuous of order $k>0$ in $y$.

\item[$B_5$:] Kolmogorov's $\varepsilon$-entropy fulfills 
\begin{align*}
\sum_{n=1}^{\infty}n^{k/2}\exp\Bigl((1-\beta)\Psi_{S_{\mathcal{E}}}(\log n/n)\Bigr)< \infty
\end{align*}
for a $\beta>1$.
\end{itemize}
Then we have the following result for the estimator of the distribution function.
\begin{theo}\label{GLTVF}
Under assumptions $B_0$ - $B_5$ we have
\begin{align*}
\underset{\chi \in S_{\mathcal{E}}}{\sup}\underset{y \in S_{\mathbb{R}}}{\sup}|\widehat{F}_{Y,B}^{\mathcal{X}}(\chi,y)-F_Y^{\mathcal{X}}(\chi,y)|=O(h_0^2)+O_{a.co.}\Bigl(\frac{\Psi_{S_{\mathcal E}}({\log n/n})}{n\ell(h_0)}\Bigr)^{1/2}.
\end{align*}
\end{theo}
The final result in this section contains the uniform rates of the bias-reduced conditional density estimator, and needs the following assumptions.
\begin{itemize}
\item[$B_4'$:] $\frac{\partial^2}{\partial y^2}f_Y^{\mathcal X}(\chi,y)$ exists for all $\chi\in S_{\mathcal E}$ and is uniformly bounded.

\item[$B_5'$:] Kolmogorov's $\varepsilon$-entropy fulfills
\begin{align*}
\sum_{n=1}^{\infty}\frac{n}{b^{3/2}}\exp((1-\beta)\Psi_{S_{\mathcal{E}}}(\log n/n))<\infty
\end{align*}
for some $\beta>1$.

\end{itemize}
\begin{theo}\label{GLTD}
Under assumptions $B_0$-$B_3$, $B_4'$ and $B_5'$, if $K_0$ is bounded, we have 
\begin{align*}
\sup_{\chi\in S_{\mathcal{E}}}\sup_{y\in S_{\mathbb{R}}}|\widehat{f}_{Y,B}^{\mathcal{X}}(\chi,y)-f_Y^{\mathcal{X}}(\chi,y)|=O(h_0^2)+O(b^2)+O_{a.co.}\Bigl(\sqrt{\frac{\Psi_{S_{\mathcal{E}}}(\log n/n)}{nb\ell(h_0)}}\Bigr).
\end{align*}

\end{theo}

\section{Finite Sample Properties}
\label{se:finite}
In this chapter we show some simulation results exemplarily for the regression model. We choose
\begin{align*}
\mathcal{X}(t)=\sin(\omega t)+t(a+2\pi)+b,~~~t\in(-1,1),
\end{align*}
with $a$ and $b$ independent identically uniformly distributed on $(0,1)$, while $\omega$ is uniformly distributed on $(0,2\pi)$ and independent of $(a,b)$. The response is
\begin{align*}
Y=r(\mathcal{X})+\varepsilon
\end{align*}
with $\varepsilon\sim\mathcal{N}(0,2)$. Let the regression function be defined as
\begin{align*}
r(\mathcal{X})=\int_{-1}^1|\mathcal{X}'(t)|(1-\cos(\pi t))dt.
\end{align*}
For computing the estimator we use a quadratic kernel $K$ and the $L_2$-metric. The fixed value $\chi$ is chosen at the start of the simulation and remains fixed for all iterations with $a=0.3064023$, $b=0.3744585$, $\omega=3.826435$ and $\varepsilon=-0.1200122$. The bias-reduced estimator is built from $B=21$ different bandwidths centered around the one for the classical estimator. We discuss different strategies for choosing the bandwidths in the sequel. In all situations presented below 500 simulation runs are performed. In Tables \ref{Tab:1}, \ref{Tab:2} and \ref{Tab:3} we compare squared bias, MSE and variance of the nonbias-reduced estimator, denoted by $\hat r$, and the bias-reduced estimator, denoted by $\hat r_B$, for different sample sizes $n$ and choices of bandwidths. The sample size $n$ is given in the left column. It can be seen, that the bias, especially of the unconstrained estimator, is considerably large even for larger sample sizes. The bias reduced estimator is able to cure this to some extend by preserving the magnitude of variance. This results in a considerably smaller MSE for all cases compared to the pilot estimator. We also see that the bandwidth choice is challenging in nonparametric functional
\begin{table}[H]
\begin{tabular}{|c|cc|cc|cc|}
\hline 
\multicolumn{1}{|c|}{}&\multicolumn{2}{|c|}{Squared Bias}&\multicolumn{2}{|c|}{MSE}&\multicolumn{2}{|c|}{Variance}\\~&
$\widehat{r}$& $\widehat{r}_B$&$\widehat{r}$&$\widehat{r}_B$&$\widehat{r}$&$\widehat{r}_B$\\\hline $sw=0.005$&9.184021&3.233732&9.207733
&3.27116&0.02371245&0.03742775\\\hline
$sw=0.01$&9.208351&3.155444&9.231983&3.193726&0.02363151&0.03828278
\\\hline
$sw=0.02$&9.212692&2.755409&9.236515&2.79365&0.02382332&0.03824129\\\hline
\end{tabular}
\caption{Squared Bias, MSE and Variance of $\hat r$ with bandwidth $h=1.2$ and of $\hat r_B$ using $B=21$ bandwidths centered around $h$ with stepwidth $sw$ given in the left column.}
\label{Tab:4}
~
\medskip
\\\begin{tabular}{|c|cc|cc|cc|}
\hline 
\multicolumn{1}{|c|}{}&\multicolumn{2}{|c|}{Squared Bias}&\multicolumn{2}{|c|}{MSE}&\multicolumn{2}{|c|}{Variance}\\~&
$\widehat{r}$& $\widehat{r}_B$&$\widehat{r}$&$\widehat{r}_B$&$\widehat{r}$&$\widehat{r}_B$\\\hline $sw=0.005$&7.380457&0.06282056&7.402475&0.1381973&0.02201808&0.07537675\\\hline
$sw=0.01$&7.262152&0.03055685&7.28609&0.1065919&0.02393775&0.07603501\\\hline
$sw=0.02$&7.306792&0.0009018255&7.330331&0.07009564&0.02353859&0.06919381\\\hline
\end{tabular}
\caption{Squared Bias, MSE and Variance of $\hat r$ with bandwidth $h=1$ and of $\hat r_B$ using $B=21$ bandwidths centered around $h$ with stepwidth $sw$ given in the left column.}
\label{Tab:5}
\end{table}

regression and that applying the bias reduction technique makes the procedure more robust against bandwidth choice. But since it is an asymptotic result, we would like to point out, that the reduction is not always visible in finite samples. This is due to the dominant effect of additional constants for the bias reduced estimator. It especially has an influence for smaller bandwidth where the bias of the classical estimator is already comparable small. For example, if $h=3.71n^{-1/3}$, the bias reduction is not visible in the above setting for $n=100,200$. For $n=500$, we have a squared bias of $0.2071$ for $\hat r$ and  $0.1184$ for $\hat r_B$. 
The results, given in Tables \ref{Tab:4} and \ref{Tab:5}, provide information about whether to choose the set of bandwidths of $\widehat{r}_B$ close to or wider around the bandwidth of $\widehat{r}$. For that we fix $n=500$.
We see, that for a starting bandwidth $h$, a wider range of bandwidths yields a better reduction of the bias. To enquire if the reduction comes from the choice of a larger stepwidth or the extension of the whole bandwidth interval we next vary the number $B$ of bandwidths used by first keeping the stepwidth constant while enlarging the bandwidth interval and by secondly keeping the bandwidth interval fixed while refining the grid of bandwidths.
\begin{table}[H]
\begin{tabular}{|c|cc|cc|cc|}
\hline 
\multicolumn{1}{|c|}{}&\multicolumn{2}{|c|}{Squared Bias}&\multicolumn{2}{|c|}{MSE}&\multicolumn{2}{|c|}{Variance}\\~&
$\widehat{r}$& $\widehat{r}_B$&$\widehat{r}$&$\widehat{r}_B$&$\widehat{r}$&$\widehat{r}_B$\\\hline
$B=11$&9.202894&3.251978&9.227184&3.289342&0.024289&0.037363\\\hline
$B=15$&9.205274&3.170393&9.226343&3.204100&0.021068&0.033706\\\hline
$B=21$&9.169364&3.120448&9.190970&3.155109&0.021605&0.034660\\\hline 
$B=41$&9.225849&2.788351&9.251343&2.832366&0.025494&0.044015\\\hline
\end{tabular}
\caption{Squared Bias, MSE and Variance of $\hat r$ with bandwidth $h=1.2$ and of $\hat r_B$ with a stepwidth of $0.01$ between two bandwidths. The left column shows how many bandwidths are used for $\hat r_B$.}
\label{Tab:6}
~
\medskip
\\\begin{tabular}{|c|cc|cc|cc|}
\hline 
\multicolumn{1}{|c|}{}&\multicolumn{2}{|c|}{Squared Bias}&\multicolumn{2}{|c|}{MSE}&\multicolumn{2}{|c|}{Variance}\\~&
$\widehat{r}$& $\widehat{r}_B$&$\widehat{r}$&$\widehat{r}_B$&$\widehat{r}$&$\widehat{r}_B$\\\hline
$B=11$&7.319643&0.056945&7.341007&0.125366&0.021364&0.068421\\\hline
$B=15$&7.344617&0.052797&7.370343&0.136344&0.025726&0.083546\\\hline
$B=21$&7.383378&0.046844&7.408037&0.114998&0.024658&0.068154\\\hline 
$B=41$&7.342917&0.002373&7.365203&0.071804&0.022286&0.069431\\\hline
\end{tabular}
\caption{Squared Bias, MSE and Variance of $\hat r$ with bandwidth $h=1$ and of $\hat r_B$ with a stepwidth of $0.01$ between two bandwidths. The left column shows how many bandwidths are used for $\hat r_B$.}
\label{Tab:7}
\end{table}
\begin{table}[H]
\begin{tabular}{|c|cc|cc|cc|}
\hline 
\multicolumn{1}{|c|}{}&\multicolumn{2}{|c|}{Squared Bias}&\multicolumn{2}{|c|}{MSE}&\multicolumn{2}{|c|}{Variance}\\~&
$\widehat{r}$& $\widehat{r}_B$&$\widehat{r}$&$\widehat{r}_B$&$\widehat{r}$&$\widehat{r}_B$\\\hline
$B=11$&9.158861&3.097066&9.182459&3.135653&0.023598&0.038587\\\hline
$B=15$&9.140199&3.089085&9.160947&3.124537&0.020747&0.035451\\\hline
$B=21$&9.214916&3.139630&9.238012&3.175078&0.023096&0.035448\\\hline 
$B=41$&9.135028&3.079761&9.157966&3.114933&0.022938&0.035172\\\hline
\end{tabular}
\caption{Squared Bias, MSE and Variance of $\hat r$ with bandwidth $h=1.2$ and of $\hat r_B$ with a stepwidth of $0.1/((B-1)/2$ between two bandwidths. The left column shows how many bandwidths are used for $\hat r_B$.}
\label{Tab:8}
~
\medskip
\\\begin{tabular}{|c|cc|cc|cc|}
\hline 
\multicolumn{1}{|c|}{}&\multicolumn{2}{|c|}{Squared Bias}&\multicolumn{2}{|c|}{MSE}&\multicolumn{2}{|c|}{Variance}\\~&
$\widehat{r}$& $\widehat{r}_B$&$\widehat{r}$&$\widehat{r}_B$&$\widehat{r}$&$\widehat{r}_B$\\\hline
$B=11$&7.338061&0.036448&7.363608&0.105866&0.025546&0.069418\\\hline
$B=15$&7.349646&0.041137&7.374302&0.120624&0.024656&0.079487\\\hline
$B=21$&7.340968&0.035893&7.365058&0.116017&0.024089&0.080123\\\hline 
$B=41$&7.383122&0.051347&7.408013&0.128210&0.024891&0.076862\\\hline
\end{tabular}
\caption{Squared Bias, MSE and Variance of $\hat r$ with bandwidth $h=1$ and of $\hat r_B$ with a stepwidth of $0.1/((B-1)/2$ between two bandwidths. The left column shows how many bandwidths are used for $\hat r_B$.}
\label{Tab:9}
\end{table}

We take different amounts of bandwidths $B\in\{11,15,21,41\}$, set $n=500$ and choose $0.01$ as the distance between two bandwidths. Results are shown in Tables \ref{Tab:6} and \ref{Tab:7}. It seems that the larger the choice of $B$ is, the better a bias reduction can be observed. Of course, a larger choice of $B$ increases the amount of calculations. {In Tables \ref{Tab:8} and \ref{Tab:9} we vary $B$, but take a constant length for the bandwidth interval, which means that also the stepwidths between the bandwidths get smaller with increasing $B$. The results show no improvement of the bias reduction. This is in line with the theoretical results in Remark \ref{rem:Bchoice}. Because of that we conclude that the length of the interval from which the bandwidths are chosen is the key for improving the biasreduction. This also holds for a quite wide interval of length 1. For a simulation with $h=2,~n=100,~B=21$ and $sw=0.05$ we get $11.28986$ as a squared bias for $\hat r$ and $9.242853$ as a squared bias for $\widehat{r}_B$, which is slightly better than the comparable result in Table \ref{Tab:1}.} We end this section by comparing the equidistant choice of bandwidths with the choice based on ideas from optimal design. To this end, in each simulation run, we take $n=500$ observations. The bandwidth $h=1$ is used for the classical estimator $\hat r$, for the bias reduced estimator with equidistant design we chose $B=20$ bandwidths $h_i=0.9+0.2\frac i{20}$, while for the bias reduced estimator with optimal design $\widehat r_{B,opt}$ we choose $h_i\in\lbrace (0.9,...,0.91),(1.09,...,1.1)\rbrace$. From the results in Table \ref{tab:design_comp} we see, that the classical estimator has the smallest variance but a considerably larger squared bias. The bias reduction for both bias reduced estimators is nearly the same, while the variance of $\widehat r_{B,opt}$ is even a little bit smaller.

\section{Concluding remarks}\label{se:out}
In this paper we proposed a general method for bias reduction of regularised estimators. At the example of nonparametric functional data models we emphasised the advantages. It turned out, that the asymptotic behavior of the classical estimators transfer to the bias reduced ones, which is of advantage in further inference like confidence regions or hypothesis testing. While the bias is reliably reduced, the finite sample variance might be slightly higher
\begin{table}[H]
\begin{tabular}{|ccc|ccc|ccc|}
\hline 
\multicolumn{3}{|c|}{Squared Bias}&\multicolumn{3}{|c|}{MSE}&\multicolumn{3}{|c|}{Variance}\\
$\widehat{r}$& $\widehat{r}_B$& $\widehat{r}_{B,opt}$&$\widehat{r}$&$\widehat{r}_B$&$\widehat{r}_{B,opt}$&$\widehat{r}$&$\widehat{r}_B$&$\widehat{r}_{B,2}$\\\hline
6.4941&0.2511&0.2519&6.5241&0.3590&0.3581& 0.0299&0.1079&0.1061\\\hline
\end{tabular}
\caption{\label{tab:design_comp} Comparison of squared biases and variances of the estimators $\hat r$, $\hat r_B$ and $\hat r_{B,opt}$.}
\end{table}
for the bias reduced estimator, but still of the same asymptotic order as before. To this end we proposed to make use of results from optimal design of experiments to choose the set of regularization parameters which yields the bias reduced estimator with a variance as small as possible. While the asymptotic results are well established, the optimal choice of the number of regularisation parameters, as well as the regularisation parameters itself, might be improved in future work to get even better finite sample results.

\paragraph{Acknoledgements:} The authors would like to thank Holger Dette and Stanislav Volgushev for helpful discussions at an early stage of this project.

\bibliography{Literatur}

\appendix

\section{Appendix: Proofs}
\label{se:app}
In the following subsections we show the proofs for the main results from the above sections 
\subsection{Proof of Theorem 3.3 and Corollary 3.4}

Because of linearity we have
\begin{align*}
E[\widehat{m}_{\Phi,B}(\chi)]=\sum_{i=1}^B g_{i,j}E[\widehat{m}_{\Phi,\textbf{h}_i}(\chi)]
\end{align*}
as the $g_{i,j}$, $i=1,\ldots,B$, $j=1,2$ are constants. For the regression case $E[\widehat{m}_{\Phi,\textbf{h}_i}(\chi)]$ has been computed in \cite{Ferraty2007} and can be derived along the lines for the conditional distribution. By a second order Taylor expansion instead of a first order one, we follow the same steps to get \begin{align*}
E[\widehat{m}_{\Phi,\textbf{h}_i}(\chi)]-m_{\Phi}(\chi)=&h_i\varphi_{\Phi,\chi}'(0)\int tK(t)dP^{||\mathcal{X}-\chi||/h_i}(t)+\frac12h_i^2\varphi_{\Phi,\chi}''(0)\int t^2K(t)dP^{||\mathcal{X}-\chi||/h_i}(t))\\&+R_n(\chi,\Phi)+o(h_i^2).
\end{align*}
With similar computations as in \cite{Ferraty2007}, using $t^2K(t)=\int_0^t s^2K'(s)ds+\int_0^t 2sK(s)ds$, we get
\begin{align*}
\int t^2K(t)dP^{||\mathcal{X}-\chi||/h_i}(t)\rightarrow \int_0^1 (s^2K'(s)+2sK(s))\tau_{0,\chi}(s) ds=:M_{3,\chi}.
\end{align*}
Therefore we have
\begin{align*}
E[\widehat{m}_{\Phi,B}(\chi)]=&\sum_{i=1}^Bg_{i,1}m_{\Phi}(\chi)+\sum_{i=1}^Bg_{i,1}h_i\varphi_{\Phi,\chi}'(0)\frac{M_{0,\chi}}{M_{1,\chi}}+\frac12\sum_{i=1}^Bg_{i,1}h_i^2\varphi_{\Phi,\chi}''(0)\frac{M_{3,\chi}}{M_{1,\chi}}\\&+\sum_{i=1}^Bg_{i,1}(R_n(\chi,\Phi)+o(h_i^2)).
\end{align*}
Remark that $\sum_{i=1}^Bg_{i,1}=1$ and $\sum_{i=1}^Bg_{i,1}h_i=0$ as in (\ref{properties}) as well as
\begin{align*}
\sum_{i=1}^Bg_{i,1}h_i^2&=\sum_{i=1}^B\frac{C_i^2h_0^2(\sum_{k=1}^BC_k^2h_0^2-C_ih_0\sum_{k=1}^BC_kh_0)}{B\sum_{k=1}^BC_k^2h_0^2-(\sum_{k=1}^BC_kh_0)^2}
\\&=h_0^2\frac{(\sum_{k=1}^BC_k^2)^2-\sum_{k=1}^BC_k^3\sum_{k=1}^BC_k}{B\sum_{k=1}^BC_k^2-(\sum_{k=1}^BC_k)^2}>0.
\end{align*}
We get
\begin{align*}
E[\widehat{m}_{\Phi,B}(\chi)]=m_{\Phi}(\chi)+\frac12\frac{(\sum_{k=1}^BC_k^2)^2-\sum_{k=1}^BC_k^3\sum_{k=1}^BC_k}{B\sum_{k=1}^BC_k^2-(\sum_{k=1}^BC_k)^2}h_0^2\varphi_{\Phi,\chi}''(0)\frac{M_{3,\chi}}{M_{1,\chi}}+R_n(\chi,\Phi)+o(h_0^2).
\end{align*} 
which proves Theorem \ref{TheoBias}. For the regression and conditional distribution function $R_n(\chi,\Phi)=0$ which directily proves the assertion. For the conditional density function $\widehat{f}_{Y,B}^{\mathcal{X}}(\chi,y)$ the term $R_n(\chi,\Phi)$ has to be calculated explicitly. There is by the law of total expactation
\begin{align*}
&E[(b^{-1}K_0(b^{-1}(y-Y))-f_Y^{\mathcal{X}}(\chi,y))K(h^{-1}||\mathcal{X}-\chi||)]\\
=&E[\int_{\mathbb{R}}K_0(v)(f_Y^{\mathcal{X}}(\mathcal{X},y-vb)-f_Y^{\mathcal{X}}(\chi,y))dvK(h^{-1}||\mathcal{X}-\chi||)].
\end{align*}
With a Taylor expansion of second, respectively fourth order of $f_Y^{\mathcal{X}}(\mathcal{X},y-vb)$ in $y$ we have
\begin{align*}
f_Y^{\mathcal{X}}(\mathcal{X},y-vb)= f_Y^{\mathcal{X}}(\mathcal{X},y)-f_Y^{'\mathcal{X}}(\mathcal{X},y)vb+\frac12f_Y^{''\mathcal{X}}(\mathcal{X},y)v^2b^2+o(b^2).
\end{align*}
respectively
\begin{align*}
f_Y^{\mathcal{X}}(\mathcal{X},y-vb)= &f_Y^{\mathcal{X}}(\mathcal{X},y)-f_Y^{'\mathcal{X}}(\mathcal{X},y)vb+\frac12f_Y^{''\mathcal{X}}(\mathcal{X},y)v^2b^2
\\&-\frac{1}{6}f_Y^{'''\mathcal{X}}(\mathcal{X},y)v^3b^3+\frac{1}{24}f_Y^{(4),\mathcal{X}}(\mathcal{X},y)v^4b^4+o(b^4).
\end{align*}
Using symmetry properties of $v^kK_0(v)$ as well as linearity of the expectation this yields
\begin{align*}
E[\widehat{f}_{Y}^{\mathcal{X}}(\chi,y)]=&f_Y^{\mathcal{X}}(\chi,y)+h\varphi'_{dens,\chi}(0)\frac{M_{0,\chi}}{M_{1,\chi}}+\frac12h^2\varphi_{dens,\chi}^{''}(0)\frac{M_{3,\chi}}{M_{1,\chi}}+o(h^2)
\\&+\frac12b^2\int_{\mathbb{R}}v^2K_0(v)dvf_Y^{''\chi}(y)+\frac{1}{24}b^4\int_{\mathbb{R}}v^4K_0(v)dvf_Y^{(4),\chi}(y)+o(b^4)
\end{align*}
respectively
\begin{align*}
E[\widehat{f}_{Y,B,2}^{\mathcal{X}}(\chi,y)]=(\sum_{i=1}^Bg_{i,2})f_Y^{\mathcal{X}}(\chi,y)+\frac12(\sum_{i=1}^Bg_{i,2}h_i)\varphi'_{dens,\chi}(0)\frac{M_{0,\chi}}{M_{1,\chi}}+(\sum_{i=1}^Bg_{i,2}h_i^2)\varphi_{dens,\chi}^{''}(0)\frac{M_{3,\chi}}{M_{1,\chi}}
\\+(\sum_{i=1}^Bg_{i,2}b_i^2)\frac12\int_{\mathbb{R}}v^2K_0(v)dvf_Y^{''\chi}(y)+(\sum_{i=1}^Bg_{i,2}b_i^4)\frac{1}{24}\int_{\mathbb{R}}v^4K_0(v)dvf_Y^{(4),\chi}(y)
\\+o(\sum_{i=1}^Bg_{i,2}h_i^2)+o(\sum_{i=1}^Bg_{i,2}b_i^4).
\end{align*}
We rewrite
\begin{align*}
\sum_{i=1}^Bg_{i,2}h_i^2=h_0^2\sum_{i=1}^Bw_{i,2}C_i^2 \text{ and } \sum_{i=1}^Bg_{i,2}b_i^4=b_0^4\sum_{i=1}^Bw_{i,2}C_i^{'4}.
\end{align*}
which concludes the proof. $\hfill \Box$

\subsection{Proof of Theorem 3.7}
We rewrite
\begin{align*}
\widehat{m}_{\Phi,B}(\chi)&=\sum_{i=1}^Bg_{i,j}\sum_{k=1}^n\frac{K(h_i^{-1}||\mathcal{X}_k-\chi||)\Phi_n(Y_k)}{\sum_{l=1}^nK(h_i^{-1}||\mathcal{X}_l-\chi||)}
\\&=\sum_{k=1}^nG_k(\mathcal X_1,\ldots,\mathcal X_n)\Phi_n(Y_k)
\end{align*}
with
\begin{align*}
G_k(\mathcal X_1,\ldots,\mathcal X_n)=\sum_{i=1}^Bg_{i,j}\frac{\frac{1}{nL_{\chi}(h_i)}K(h_i^{-1}||\mathcal{X}_k-\chi||)}{\frac{1}{nL_{\chi}(h_i)}\sum_{l=1}^nK(h_i^{-1}||\mathcal{X}_l-\chi||)}
\end{align*}
and in the sequel $K_{i,k}:=K(h_i^{-1}||\mathcal{X}_k-\chi||)$. Because of the law of total variance we can write
\begin{align*}
Var(\widehat{m}_{\Phi,B}(\chi))=Var(E[\widehat{m}_{\Phi,B}(\chi)|\mathcal{X}_1,...,\mathcal{X}_n])+E[Var(\widehat{m}_{\Phi,B}(\chi)|\mathcal{X}_1,...,\mathcal{X}_n)].
\end{align*}
By eventually neglecting remainder terms, the first term is estimated as
\begin{align*}
Var(E&[\widehat{m}_{\Phi,B}(\chi)|\mathcal{X}_1,...,\mathcal{X}_n])=Var(\sum_{k=1}^nG_k(\mathcal{X}_1,...,\mathcal{X}_n)E[\Phi_n(Y_k)|\mathcal{X}_k])
\\\sim &\sum_{k=1}^nVar(G_k(\mathcal{X}_1,...,\mathcal{X}_n)m_{\Phi}(\mathcal{X}_k))\\
=&\sum_{k=1}^nVar(G_k(\mathcal{X}_1,...,\mathcal{X}_n)(m_{\Phi}(\mathcal{X}_k)-m_{\Phi}(\chi)))\\
\leq& \frac{1}{n}E[\Bigl(\sum_{i=1}^Bg_{i,j}\frac{\frac{1}{L_{\chi}(h_i)}K_{i,1}}{\frac{1}{nL_{\chi}(h_i)}\sum_{l=1}^nK_{i,l}}\Bigr)^2(m_{\Phi}(\mathcal{X}_1)-m_{\Phi}(\chi))^2]\\
=&\frac{1}{nL_{\chi}(h_0)}\sum_{i_1,i_2=1}^Bg_{i_1,j}g_{i_2,j}\frac{1}{L_{\chi}(h_{i_2})}\tau_{h_{i_1},\chi}(1/C_{i_1})\\&\times E\left[(m_{\Phi}(\mathcal{X}_1)-m_{\Phi}(\chi))^2\frac{K_{i_1,1}K_{i_2,1}}{\frac{1}{n^2L_{\chi}(h_{i_1})L_{\chi}(h_{i_2})}\sum_{l=1}^nK_{i,l}\sum_{k=1}^nK_{j,k}}\right]\\
\sim & \frac1{nL_{\chi}(h_0)}\frac{1}{M_{1,\chi}^2}E[\overline{\varphi}(||\mathcal{X}_1-\chi||)K_{i_1,1}K_{i_2,1}] = o\left(\frac1{nL_{\chi}(h_0)}\right)
\end{align*}
by continuity of the function
\[\bar\varphi(s)=E[(m_{\Phi}(\mathcal X_1)-m_{\Phi}(\chi))^2|\, ||\mathcal X_1-\chi||=s].\]
Regarding the second term we have
\begin{align*}
E[Var&(\widehat{m}_{\Phi,B}(\chi)|\mathcal{X}_1,...,\mathcal{X}_n)]=E[\sum_{k=1}^nVar(G_k(\mathcal{X}_1,...,\mathcal{X}_n)\Phi_n(\mathcal X_k)|\mathcal{X}_1,...,\mathcal{X}_n)]
\\=&\sum_{k=1}^nE[G_k^2(\mathcal{X}_1,...,\mathcal{X}_n)Var(\Phi_n(\mathcal X_k)|\mathcal{X}_1,...,\mathcal{X}_n)]
\\=&\sum_{k=1}^nE[G_k^2(\mathcal{X}_1,...,\mathcal{X}_n)D_{\Phi,n}(\mathcal X_k)]\\
\sim&\frac{1}{n}\sum_{i_1,i_2=1}^Bg_{i_1,1}g_{i_2,1}\frac{1}{L_{\chi}(h_{i_1})L_{\chi}(h_{i_2})}\frac{E[K_{i_1,1}K_{i_2,1}]}{E[\frac{1}{n^2L_{\chi}(h_{i_1})L_{\chi}(h_{i_2})}\sum_{l=1}^nK_{i_1,l}\sum_{k=1}^nK_{i_2,k}]}D_{\Phi,n}(\chi)
\end{align*}
by continuity, see e.g.~Lemma 3 in \cite{Ferraty2007}.

Because of the i.i.d. assumption, the denominators in the weights can be written as
\begin{align*}
E[&\frac{1}{n^2L_{\chi}(h_{i_1})L_{\chi}(h_{i_2})}\sum_{l=1}^nK_{i_1,l}\sum_{k=1}^nK_{i_2,k}]\\
&=\frac{1}{n^2}((n^2-n)E\left[\frac1{L_{\chi}(h_{i_1})}K_{i_1,1}\right]E\left[\frac1{L_{\chi}(h_{i_2})}K_{i_2,1}\right]+nE\left[\frac1{L_{\chi}(h_{i_1})L_{\chi}(h_{i_2})}K_{i_1,1}K_{i_2,1}\right])\\
&\to M_{1,\chi}^2
\end{align*}
for $n\to\infty$. 

To compute the expectation of kernels in the numerator, we follow the same steps as in \cite{Ferraty2007}.
\begin{align*}
E[K_{i_1,1}K_{i_2,1}]=\int K(t/h_{i_1})K(t/h_{i_2})dP^{||\mathcal{X}-\chi||}(t)=\int K(t)K\left(t\frac{h_{i_2}}{h_{i_1}}\right)dP^{||\mathcal{X}-\chi||/h_{i_2}}(t).
\end{align*}
Using $C_{i_2}h_0=h_{i_2}$ and Fubini we have
\begin{align*}
&\int K(t)K(t\frac{h_{i_1}}{h_{i_1}})dP^{||\mathcal{X}-\chi||/h_{i_2}}(t)
\\=&K(1)K(\frac{C_{i_2}}{C_{i_1}})L_{\chi}(h_{i_2})-\int_0^1\left(\int_t^1(K(s)K'\left(s\frac{C_{i_2}}{C_{i_1}}\right)L_\chi(h_{i_2}s)ds\right)dP^{||\mathcal{X}-\chi||/h_{i_2}}
\\=&K(1)K(\frac{C_{i_2}}{C_{i_1}})L_{\chi}(h_{i_2})-\int_0^1(K(s)K'\left(s\frac{C_{i_2}}{C_{i_1}}\right)L_\chi(h_{i_2}s)ds.
\end{align*}
With 
\begin{align}\label{jauf0}
\frac{1}{L_\chi(h_{i})}=\frac{1}{L_\chi(h_0)}\frac{L_\chi(h_0)}{L_\chi(h_{i})}=\frac{1}{L_\chi(h_0)}\frac{L_\chi(\frac{1}{C_{i}}h_{i})}{L_\chi(h_{i})}=\frac{1}{L_\chi(h_0)}\tau_{h_{i},\chi}(\frac{1}{C_{i}}),
\end{align}
and $\tau_{h_i,\chi}(1/C_i)\rightarrow\tau_{0,\chi}(1/C_i)$ it holds that
\begin{align*}
M_{2,i_1,i_2,\chi}&=\lim_{h_{i_2}\to 0}\frac{1}{L_\chi(h_{i_2})}\left\{K(1)K(C_{i_2}/C_{i_1})L_\chi(h_{i_2})-\int_0^1K(s)K'\left(s\frac{C_{i_2}}{C_{i_1}}\right)L_\chi(h_{i_2}s)ds\right\}
\\&=K(1)K\left(\frac{C_{i_2}}{C_{i_1}}\right)-\int_0^1K(s)K'\left(s\frac{C_{i_2}}{C_{i_1}}\right)\tau_{0,\chi}(s)ds.
\end{align*}
Since the weights $g_{i,j}$, $i=1,\ldots,B$, $j=1,2$ are independent of $n$ , it follows that
\begin{align*}
Var(\hat m_{\Phi,B}(\chi))\sim \frac1{nL_{\chi}(h_0)}\frac{D_{\Phi,n}(\chi)}{M_{1,\chi}^2}\sum_{i_1,i_2=1}^Bg_{i_1,j}g_{i_2,j}\tau_{0,\chi}(1/C_{i_1})M_{2,i_1,i_2,\chi}.
\end{align*}

\subsection{Proof of Theorem 3.8}
Since we know from Theorem \ref{TheoVar} that $Var(\widehat m_{\Phi,B}(\chi))=O(1/nL_{\chi}(h_0))$, $L_\chi(h_i)=O(L_\chi(h_0))$ and $\frac1{nL_\chi(h_i)}\sum_{l=1}^nK(h_i^{-1}||\mathcal{X}_l-\chi||)\sim M_{1,\chi}$, it suffices to show
\begin{align*}
\frac1{n^4L_\chi^4(h_0)}&\sum_{k=1}^n E[|\Phi_n(Y_k)\sum_{i=1}^B g_{i,j}{K(h_i^{-1}||\mathcal{X}_k-\chi||)}-E[\Phi_n(Y_k)\sum_{i=1}^B g_{i,j}{K(h_i^{-1}||\mathcal{X}_k-\chi||)}]|^4]\\
&=o\left(\frac{1}{(nL_\chi(h_0))^2}\right).
\end{align*}
Using the Minkowski and Jensen inequalities, the law of total expectation, as well as the fact, that $E[\Phi_n^4(Y_k)|\mathcal X_k=\chi]$ is bounded, we get for some constant $0<\kappa<\infty$
\begin{align*}
\frac1{n^4L_\chi^4(h_0)}&\sum_{k=1}^n E[|\Phi_n(Y_k)\sum_{i=1}^B g_{i,j}{K(h_i^{-1}||\mathcal{X}_k-\chi||)}-E[\Phi_n(Y_k)\sum_{i=1}^B g_{i,j}{K(h_i^{-1}||\mathcal{X}_k-\chi||)}]|^4]\\
&\leq \frac\kappa{n^4L_\chi^4(h_0)}\sum_{k=1}^n\sum_{i=1}^B( |g_{i,j}|^4\left(E[{K^4(h_i^{-1}||\mathcal{X}_k-\chi||)}-\left(E[{K(h_i^{-1}||\mathcal{X}_k-\chi||)}]\right)^4\right)\\
&=\frac{B\kappa}{n^3L_\chi^4(h_0)}\left(O(L_\chi(h_0))+O(L_\chi^4(h_0))\right)=o\left(\frac{1}{(nL_\chi(h_0))^2}\right).
\end{align*}

$\hfill \Box$

\subsection{Proof of Theorems 3.9, 3.10 and 3.11}\label{ProofGKR}
Again the proofs for all three estimators are very similar. Therefore we discuss them all in this section. Let $\hat m_{\Phi,h_i}$ denote the estimator of $m_\Phi$ with bandwidth $h_i$. We then obtain the decomposition
\begin{align*}
\widehat{m}_{\Phi,B}(\chi)&-m_\Phi(\chi)=\sum_{i=1}^Bg_{i,j}\widehat{m}_{\Phi,h_i}(\chi)-m_\Phi(\chi)=\sum_{i=1}^Bg_{i,j}\frac{\widehat{e}_{\Phi,h_i}(\chi)}{\widehat{f}_{h_i}(\chi)}-m_\Phi(\chi)\\
=&\sum_{i=1}^Bg_{i,j}\frac1{\hat f_{h_i}}(\hat e_{\Phi,h_i}(\chi)-E[\hat e_{\Phi,h_i}(\chi)])+\sum_{i=1}^Bg_{i,j}(E[\hat e_{\Phi,h_i}(\chi)]-m_\Phi(\chi))\\&+\sum_{i=1}^Bg_{i,j}\frac{E[\hat e_{\Phi,h_i}(\chi)]}{\hat f_{h_i}(\chi)}(1-\hat f_{h_i}(\chi))\\
=&T_1(\chi)+T_2(\chi)+T_3(\chi).
\end{align*}
with
\begin{align*}
\widehat{f}_{h_i}(\chi)&=\frac{1}{nE[K(h_i^{-1}||\mathcal{X}-\chi||)]}\sum_{k=1}^nK(h_i^{-1}||\mathcal{X}_k-\chi||)\\
\widehat{e}_{\Phi,h_i}(\chi)&=\frac{1}{nE[K(h_i^{-1}||\mathcal{X}-\chi||)]}\sum_{k=1}^nK(h_i^{-1}||\mathcal{X}_k-\chi||)\Phi_n(Y_k).
\end{align*}
Using condition $B_1$, which yields $\ell(h_i)\sim\ell(h_0)$, and the fact, that $\sum_{i=1}^B|g_{i,j}|<\infty$, we obtain by Lemmas 8 and 11, 13 respectively 15 as well as Corollary 9 in \cite{Ferraty2010} 
\begin{align*}
\sup_{\chi\in S_{\mathcal E}}|T_3(\chi)|&=O_{a.co.}\left(\frac{\Psi_{S_{\mathcal E}}({\log n/n})}{n\ell(h_0)}\right)^{1/2}\\
\sup_{\chi\in S_{\mathcal E}}|T_1(\chi)|&=O_{a.co.}\left(\frac{\Psi_{S_{\mathcal E}}({\log n/n})}{nd_{\Phi,n}\ell(h_0)}\right)^{1/2}
\end{align*}
and for the conditional distribution function respectively the conditional density function
\begin{align*}
\sup_{\chi\in S_{\mathcal E},y\in S_{\mathbb R}}|T_1(\chi)|&=O_{a.co.}\left(\frac{\Psi_{S_{\mathcal E}}({\log n/n})}{nd_{\Phi,n}\ell(h_0)}\right)^{1/2}.\\
\end{align*}
For the bias term $T_2(\chi)$, Theorem \ref{TheoBias} yields
\begin{align*}
\sup|T_2(\chi)|&=\sup\Big|S(\chi,\Phi)h_0^2+R_n(\chi,\Phi)+O((nL_\chi(h))^{-1})+o(h_0^2)\Big|\\
&\leq \sup|S(\chi,\Phi)|h_0^2+\sup|R_n(\chi,\Phi)|+O((n\ell(h_0))^{-1})+o(h_0^2)\\
&=O(h_0^2)+\sup|R_n(\chi,\Phi)|+O((n\ell(h_0))^{-1})+o(h_0^2)
\end{align*}
where the $\sup$ is over $\{\chi\in E_{\mathcal E}\}$ respectively $\{\chi\in E_{\mathcal E},y\in S_{\mathbb R}\}$. This proves Theorems \ref{GKT}-\ref{GLTD}.\hfill $\Box$

\end{document}